\documentclass[]{interact}

\usepackage{epstopdf}
\usepackage[caption=false]{subfig}

\usepackage[numbers,sort&compress]{natbib}
\bibpunct[, ]{[}{]}{,}{n}{,}{,}
\makeatletter
\def\NAT@def@citea{\def\@citea{\NAT@separator}}
\makeatother

\theoremstyle{plain}
\newtheorem{theorem}{Theorem}[section]

\newtheorem{proposition}[theorem]{Proposition}

\theoremstyle{definition}

\theoremstyle{remark}
\newtheorem{remark}{Remark}

\begin{document}


\title{Properties of the multi-index special function $\mathcal{W}^{\left(\bar{\alpha},\bar{\nu}\right)}(z)$}

\author{
\name{R. Droghei\textsuperscript{a}\thanks{CONTACT R. Droghei. Email: riccardo.droghei@francescoseveri.org} }
\affil{\textsuperscript{a}Liceo Scientifico Francesco Severi,
Viale Europa,36, 03100 Frosinone (FR), ITALY }
}

\maketitle

\begin{abstract}
In this paper, we investigate some properties related to a multi-index special function $\mathcal{W}^{\left(\bar{\alpha},\bar{\nu}\right)}$ that arose from an eigenvalue problem for a multi-order fractional hyper-Bessel operator, involving Caputo fractional derivatives. We show that for particular values of the parameters involved in this special function $\mathcal{W}^{\left(\bar{\alpha},\bar{\nu}\right)}$, this leads to the hyper-Bessel function of Delerue. The Laplace transform of the $\mathcal{W}^{\left(\bar{\alpha},\bar{\nu}\right)}$ is discussed obtaining, in particular cases, the well-known functional relation between hyper-Bessel function and multi-index Mittag-Leffler function, or, quite simply, between classical Wright and Mittag-Leffler functions. Moreover, it is shown that the multi-index special function satisfies the recurrence relation involving fractional derivatives. In a particular case, we derive, to the best of our knowledge, a new differential recurrence relation for the Mittag-Leffler function. We also provide derivatives of the  3-parameters function $\mathcal{W}_{\alpha,\beta,\nu}$ with respect to parameters, leading to infinite power series with coefficients being quotients of digamma and gamma functions.
\end{abstract}


\section{Introduction}

Nowadays, the interest in fractional differential equations is increasing because these are becoming more adequate than those of integer order to investigate various problems in different fields of physics, engineering and economics \cite{Pod99}, \cite{FM2010}, \cite{GKMR2020}. They have indeed the fundamental characteristic to describe memory and heredity properties of many materials. Some of them have been introduced within the framework of partition theory in solving number theory problems. This is the case of the Wright function, introduced by E. M. Wright in his articles on the asymptotic partition formulae\cite{WA1}, \cite{WA2}, \cite{WA3} and \cite{Wright1933}, \cite{Wright1935a}, \cite{Wright1935b}.

Recently, many authors are dealing with {\it{multi-indices}} special functions (SF) of fractional calculus (FC) appearing in solution of differential equations and systems of fractional multi-order type (e.g. hyper-Bessel and quasi-Bessel operators) \cite{GarraPolito}, \cite{DS22}. Among them, the most general functions we just want to refer to are the {\it{Fox H-function}} and the {\it{Wright generalized hypergeometric function}} \cite{K2019}. Indeed, one gets the classical SF setting their parameters with integer values.

In the previous paper \cite{Droghei21} the author investigated a hyper-Bessel-type operator involving Caputo derivatives. Solving the eigenvalue problem associated with this fractional operator, the author introduced a function, written in series expansion, that in specific cases is possible to refer to the well-known special function of the fractional calculus. 
According to the information we have, this special function was not studied by now. But as seen, it is reduced in particular cases to some known special functions, which on their side are cases of the Bessel and hyper-Bessel functions and more generally, of the multi-index Mittag-Leffer functions.

This multi-index special function, called in the previous paper m-p generalized Wright function, plays an important role in nonlinear fractional differential equations, and in their isochronous $\omega$-modified version\cite{Droghei21},\cite{DG20}. It is also a natural generalization of the applications of the Laguerre derivatives and the Laguerre-type exponentials \cite{DR2003}, \cite{BR07}, \cite{Ricci2020}, \cite{GT21}. In this survey article, firstly, we want to examine several properties associated with the multi-index special function investigated in \cite{Droghei21}.

The outline of this work is as follows. In Section 2, we recall the definition of the multi-index function $\mathcal{W}^{\left(\bar{\alpha},\bar{\nu}\right)}$ introduced in \cite{Droghei21} and its connection with the Hyper -Bessel function. Moreover, the simpler function in the only 3-parameters case $\mathcal{W}_{\alpha,\beta,\nu}$ is described.
In Section 3 we computed the Laplace Transform of the function $\mathcal{W}^{\left(\bar{\alpha},\bar{\nu}\right)}$ and, using it, we derived some new functional relations between this function and other known special functions. The main result of this work is described in Section 4. Here we showed the recurrence relations of the function $\mathcal{W}_{\alpha,\beta,\nu}$  obtaining, we suppose, new differential recurrence relation for the Mittag-Leffler function. In Section 5 we investigated the derivatives of $\mathcal{W}_{\alpha,\beta,\nu}$ with respect to the parameters.

\section{Multi-index special function $\mathcal{W}^{\left(\bar{\alpha},\bar{\nu}\right)}(z)$}

The multi-index special function $\mathcal{W}^{\left(\bar{\alpha},\bar{\nu}\right)}(z)$ investigated in \cite{Droghei21}, is defined by series representation as a function of the complex variable $z$ and parameters $\alpha_j,\,\,j=1,..., n+1$ and $\nu_j,\,\,j=1,...,n$:

\begin{eqnarray}\label{mpgWcomp}
 &&\mathcal{W}^{\left(\bar{\alpha},\bar{\nu}\right)}(z)=\sum_{k=0}^\infty\prod_{i=1}^k\prod_{j=1}^n\frac{\Gamma(\alpha_{n+1} i+a_j)}{\Gamma(\alpha_{n+1} i+b_j)}\cdot\frac{z^k}{\Gamma(\alpha_{n+1} k+b_{n+1})}\nonumber\\. 
\end{eqnarray}

where 
\begin{eqnarray}\label{param}
\nonumber a_j=1+\sum_{m=1}^j\left(\nu_{m-1}-\alpha_{m}\right);\\
b_j=1+\sum_{m=1}^j\left(\nu_{m-1}-\alpha_{m-1}\right).
\end{eqnarray}
and the relation $a_j=b_j-\alpha_j$ with $j=1..n+1.$

The $\mathcal{W}^{\left(\bar{\alpha},\bar{\nu}\right)}(z)$ is an entire function for $\alpha_j>0,\,j=1..n+1;\,\,\nu_j\in \mathbb{C},\,j=1..n$ and $\alpha_0=\nu_0=0$.

\begin{theorem}
The multi-index special function $\mathcal{W}^{\left(\bar{\alpha},\bar{\nu}\right)}(\lambda x^{\alpha_{n+1}})$ with $\lambda\in\mathbb{R},\,\,x\geq0,\,\,\alpha_j>0,\,\,j=1,...,n+1$ and $\nu_j>0,\,j=1..n$
satisfy the following fractional differential equation involving fractional \textit{hyper-Bessel-type operator}.[see \cite{Droghei21} for the proof]

\begin{equation}\label{eigen}
\hat{D}_{nL}^{(\bar\alpha,\bar\nu)}\mathcal{W}^{\left(\bar{\alpha},\bar{\nu}\right)}(\lambda x^{\alpha_{n+1}})=\lambda\mathcal{W}^{\left(\bar{\alpha},\bar{\nu}\right)}(\lambda x^{\alpha_{n+1}});
\end{equation}

where 
\begin{equation}\label{fhb}
\hat{D}_{nL}^{(\bar\alpha,\bar\nu)}=x^{\sum_{s=1}^{n} (\alpha_s-\nu_s)}\frac{d^{\alpha_{n+1}}}{dx^{\alpha_{n+1}}}x^{\nu_n}\frac{d^{\alpha_{n}}}{dx^{\alpha_{n}}}x^{\nu_{n-1}}\frac{d^{\alpha_{n-1}}}{dx^{\alpha_{n-1}}}\cdots x^{\nu_1}\frac{d^{\alpha_{1}}}{dx^{\alpha_{1}}}.
\end{equation}
\end{theorem}


\subsection{Hyper-Bessel function as a particular case}
The hyper-Bessel function of Delerue (or a multi-index analogue of Bessel function) of order $d$ with indices $\mu_1,..., \mu_d$, introduced in 1953 by Delereu \cite{Del53} as a generalization of the Bessel function of the first type (see also \cite{Kir94}) is defined by

\begin{equation}
\mathcal{J}_{\mu_d}(z)=z^{-\frac{\mu_1+...+\mu_d}{d+1}}J_{\mu_d}((d+1) \sqrt[d+1]{z})=\sum_{k\ge0}\frac{(-1)^k z^k}{k! \prod_{j=1}^d\Gamma(k+\mu_j+1)}.
\end{equation}
Setting  $\alpha_j=1,\,\,\,j=1,...,n+1$ in the multi-index special function $\mathcal{W}^{\left(\bar{\alpha},\bar{\nu}\right)}$, we obtain the relation
\begin{equation}\label{HB}
\mathcal{W}^{\left(\bar{1},\bar{\nu}\right)}(z)=\prod_{j=1}^n\Gamma(1+a_j)\mathcal{J}_{an}(-z), .
\end{equation}
with $a_j$ defined in (\ref{param}). It is not surprising because the hyper-Bessel function satisfies the so-called hyper-Bessel differential operators of higher order, introduced by Dimovski and Kiryakova \cite{DK86}, \cite{DK87}, and obtained from (\ref{eigen}) setting all parameters $\alpha_j=1$ with $j=1..n+1$, i.e. derivatives of integer order.

\subsection{3-parameters function $\mathcal{W}_{\alpha,\beta,\nu}$}

In this section we analyse the simpler case of (\ref{mpgWcomp}) with $n=1$, $\alpha_2=\beta,\,\, \alpha_1=\alpha$ and $\nu_1=\nu$:
\begin{equation}\label{gW}
\mathcal{W}_{\alpha,\beta,\nu}(x^{\beta})=\sum_{k=0}^\infty\prod_{i=1}^k\frac{\Gamma(\beta i+1-\alpha)}{\Gamma(\beta i+1)}\frac{x^{\beta k}}{\Gamma(\beta k+1-\alpha+\nu)}.
\end{equation}

\begin{proposition}\label{lag}
Obviously, the above function (\ref{gW}) satisfies the following fractional differential equation
\begin{equation}
\hat{D}_{\alpha,\beta,\nu}f(x)=x^{\alpha-\nu}\frac{d^\beta}{dx^\beta}\left(x^\nu\frac{d^\alpha}{dx^\alpha}f(x)\right)=f(x),
\end{equation}
involving two fractional derivatives in the sense of Caputo of orders $\alpha,\beta\in\left(0,1\right)$. Where $$f(x)=\mathcal{W}_{\alpha,\beta,\nu}(x^\beta)$$.
\end{proposition}

\begin{remark}{\textbf{The Weinstein  and Bessel-Clifford operators}}
Setting $\alpha=\beta=1$ and $\nu=k,\,\,\,k\ge1$ the operator $\hat{D}_{\alpha,\beta,\nu}$ becomes $$\hat{D}_{1,1,k}=x B_k=x\left(\frac{d^2}{dx^2}+\frac{k}{x}\frac{d}{dx}\right)=x^{-k+1}\frac{d}{dx}x^{k}\frac{d}{dx}$$
where $B_k$ is the well known Weinstein operator (or Bessel operator) from the so-called Darboux-Weinstein relation \cite{Wein55}, \cite{KH95}. In \cite{Hay67} Hayek studied in details exactly the operator $\hat{D}_{1,1,k+1}$ calling its solution as Bessel-Clifford function of second order $C_\nu(x)=x^{-\frac{\nu-1}{2}}I_{\nu-1}(2\sqrt x)=\frac{1}{\Gamma(\nu+1)}_0F_1(\nu+1;-x)$, where $I_\nu(x)$ is the modified Bessel function of the first kind. Later, in \cite{Hay87} he introduced the two indices Bessel-Clifford functions of the third order modifying the hyper-Bessel function $J_{\mu,\nu}^{(2)}(x)$:
\begin{equation}\label{BC3}
C_{\mu,\nu}(x)=x^{-\frac{\mu+\nu}{3}}J_{\mu,\nu}^{(2)}(3 \sqrt[3] x)=\frac{1}{\Gamma(\mu+1)\Gamma(\nu+1)}_0F_2(\mu+1,\nu+1;-x);
\end{equation} 
satisfying the third-order Bessel-Clifford differential equation related to the operator
\begin{equation}
\hat{B}_{\mu,\nu}=x^{-\nu}\frac{d}{dx}x^{\mu-\nu+1}\frac{d}{dx}x^{\nu+1}\frac{d}{dx}.
\end{equation}

As it is simple to see, the two-parameter operator $\hat{B}_{\mu,\nu}$ is equivalent to the operator (\ref{fhb}),  $\hat{D}_{2L}^{(\{\alpha_1,\alpha_2,\alpha_3\},\{\nu_1,\nu_2\})}$ with $\alpha_1=\alpha_2=\alpha_3=1;\,\,\nu_1=\nu+1$ and $\nu_2=\mu-\nu+1$; and then the Bessel-Clifford of the third order function (\ref{BC3}) is equal to $$C_{\mu,\nu}(x)=\frac{1}{\Gamma(\nu+1)}\mathcal{W}^{\left(\{1,1,1\},\{\nu+1,\mu-\nu+1\}\right)}(x).$$

These differential operators appear very often in the PDEs of mathematical physics (especially in fluid mechanics, elasticity, and transonic flow), for instance in the generalized Bessel heat equation and other equations of generalized axially symmetric potentials (GASP) theory \cite{Wein53}. 
\end{remark}

\subsubsection{\textbf{Particular cases of $\mathcal{W}_{\alpha,\beta,\nu}$}}

For $\alpha=1$, $\beta=\lambda$ and $\nu=\mu$ the function corresponds to the Classical Wright function \\
\begin{equation}\label{Wright}
\mathcal{W}_{1,\lambda,\mu}(x^\lambda)=W_{\lambda,\mu}\left(\frac{x^\lambda}{\lambda}\right)=\sum_{k=0}^{\infty}\frac{\left(\frac{x^\lambda}{\lambda}\right)^{k}}{k!\Gamma(\lambda k+\mu)}.
\end{equation}
For $\alpha=0,\,\,\beta\rightarrow\alpha,\,\,\nu\rightarrow\beta-1$ the function corresponds to the generalized Mittag-Leffler function
\begin{equation}
\mathcal{W}_{0,\alpha,\beta-1}(z)=E_{\alpha,\beta}(z)=\sum_{k=0}^{\infty}\frac{z^{k}}{\Gamma(\alpha k+\beta)}.
\end{equation}
In case of $\alpha=\beta=\nu$ holds the relation $$\mathcal{W}_{\nu,\nu,\nu}(x^\nu)=E_{1;\nu,1}(x^\nu)$$ 
where $E_{\alpha;\nu,\gamma}(x)=\sum_{k=0}^{\infty}\frac{x^k}{\Gamma^{\alpha+1}(\nu k+\gamma)}$ is the $\alpha$\textit{-Mittag-Leffler function}.
\\

\begin{figure}[h]
\centering
\subfloat[Plot of the function $\mathcal{W}_{0,1,\nu}(x)$ for $\nu=0;0.25;0.5;0.75;1;1.25;1.5;1.75;2.$]{%
\resizebox*{7cm}{!}{\includegraphics{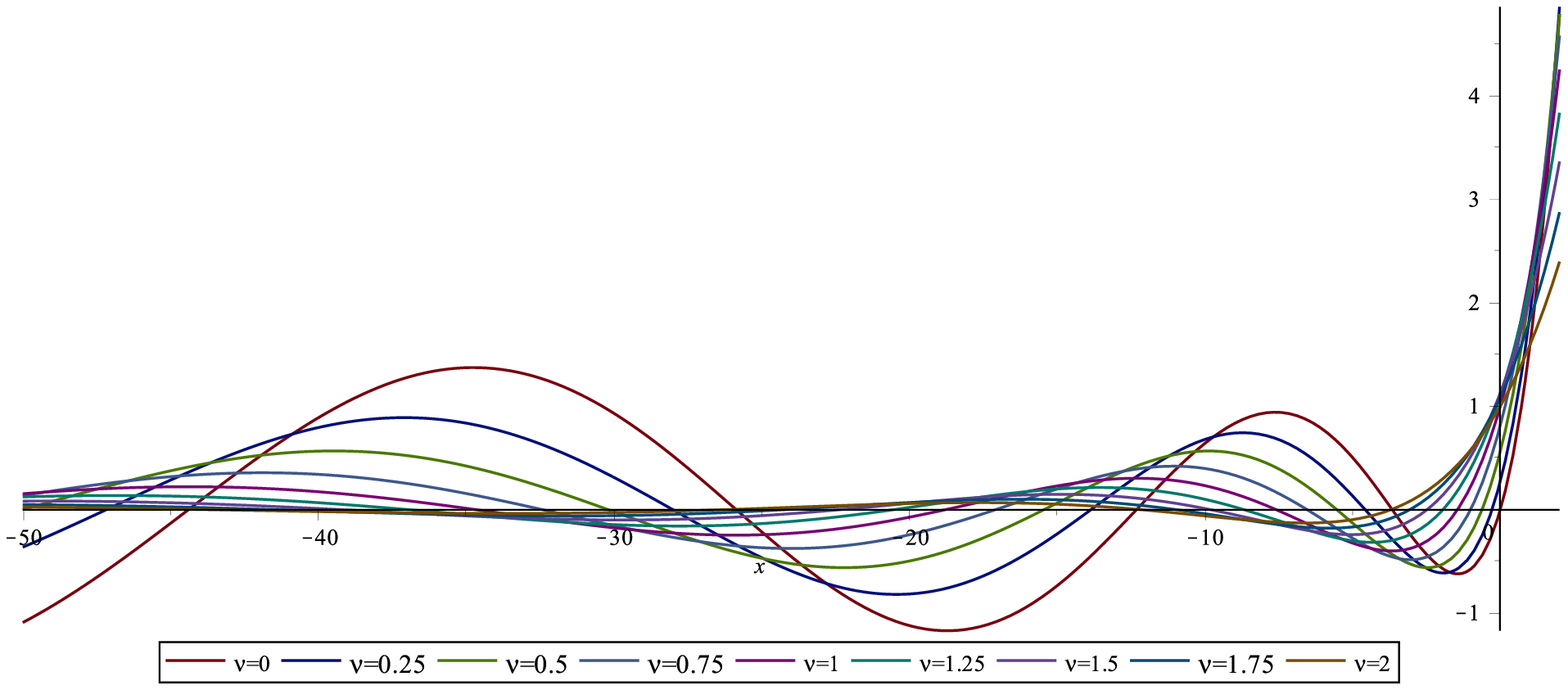}}}\hspace{5pt}
\subfloat[Plot of the function $\mathcal{W}_{1,1,\nu}(x)$ for $\nu=0;0.25;0.5;0.75;1;1.25;1.5;1.75;2.$]{%
\resizebox*{7cm}{!}{\includegraphics{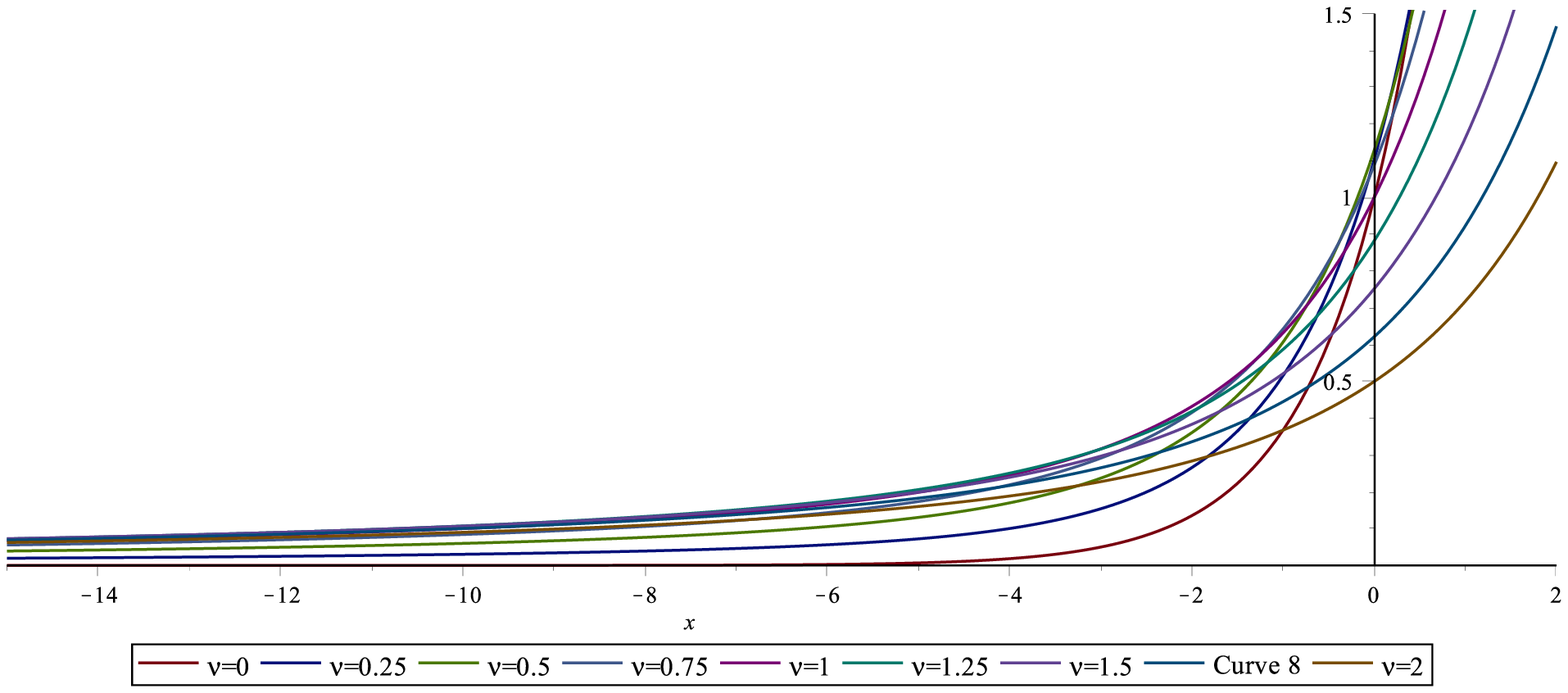}}}\hspace{5pt}
\subfloat[Plot of the function $\mathcal{W}_{\frac{1}{2},\frac{1}{2},\nu}(\sqrt x)$ for $\nu=0;0.25;0.5;0.75;1;1.25;1.5;1.75;2.$]{%
\resizebox*{7cm}{!}{\includegraphics{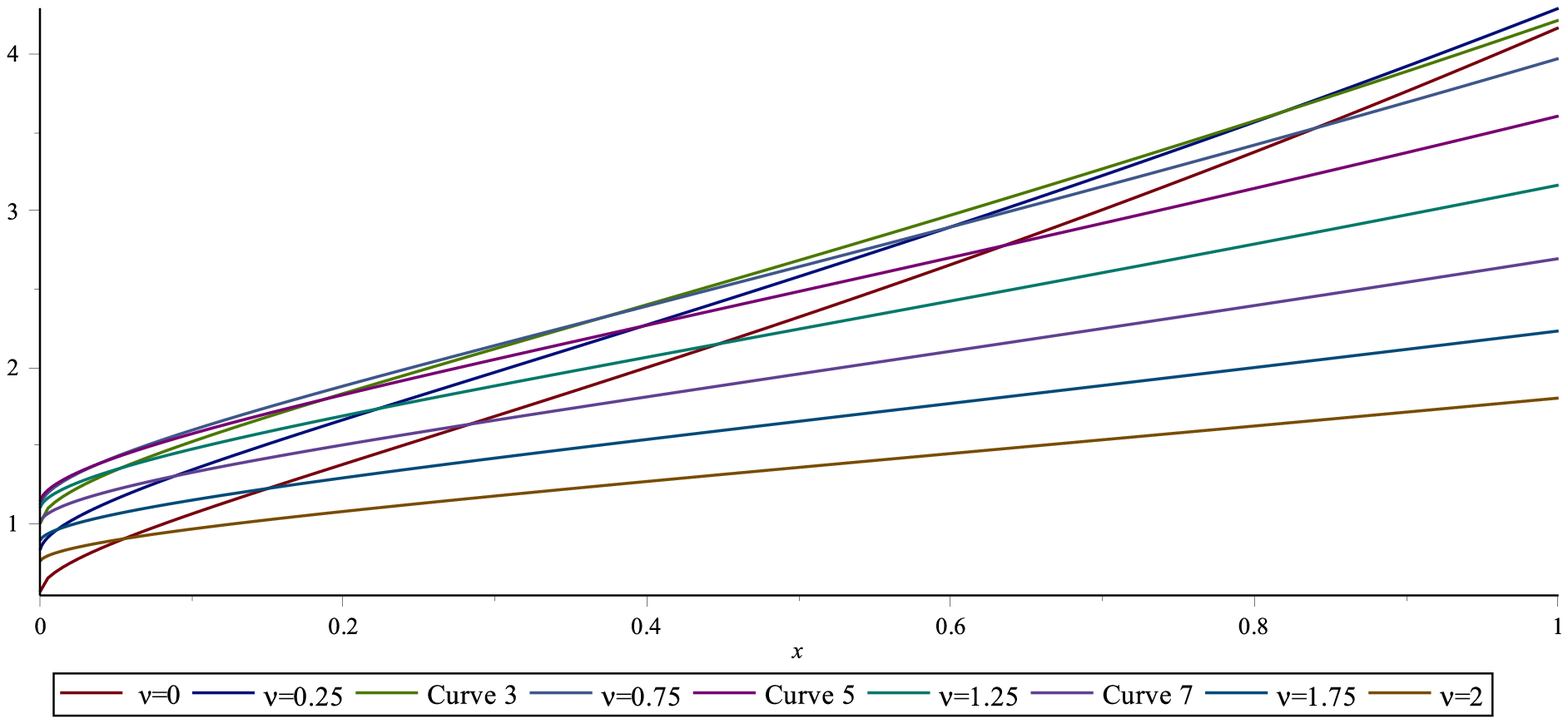}}}\hspace{5pt}
\subfloat[Plot of the function $\mathcal{W}_{\frac{1}{2},1,\nu}(x)$ for $\nu=0;0.25;0.5;0.75;1;1.25;1.5;1.75;2.$]{%
\resizebox*{7cm}{!}{\includegraphics{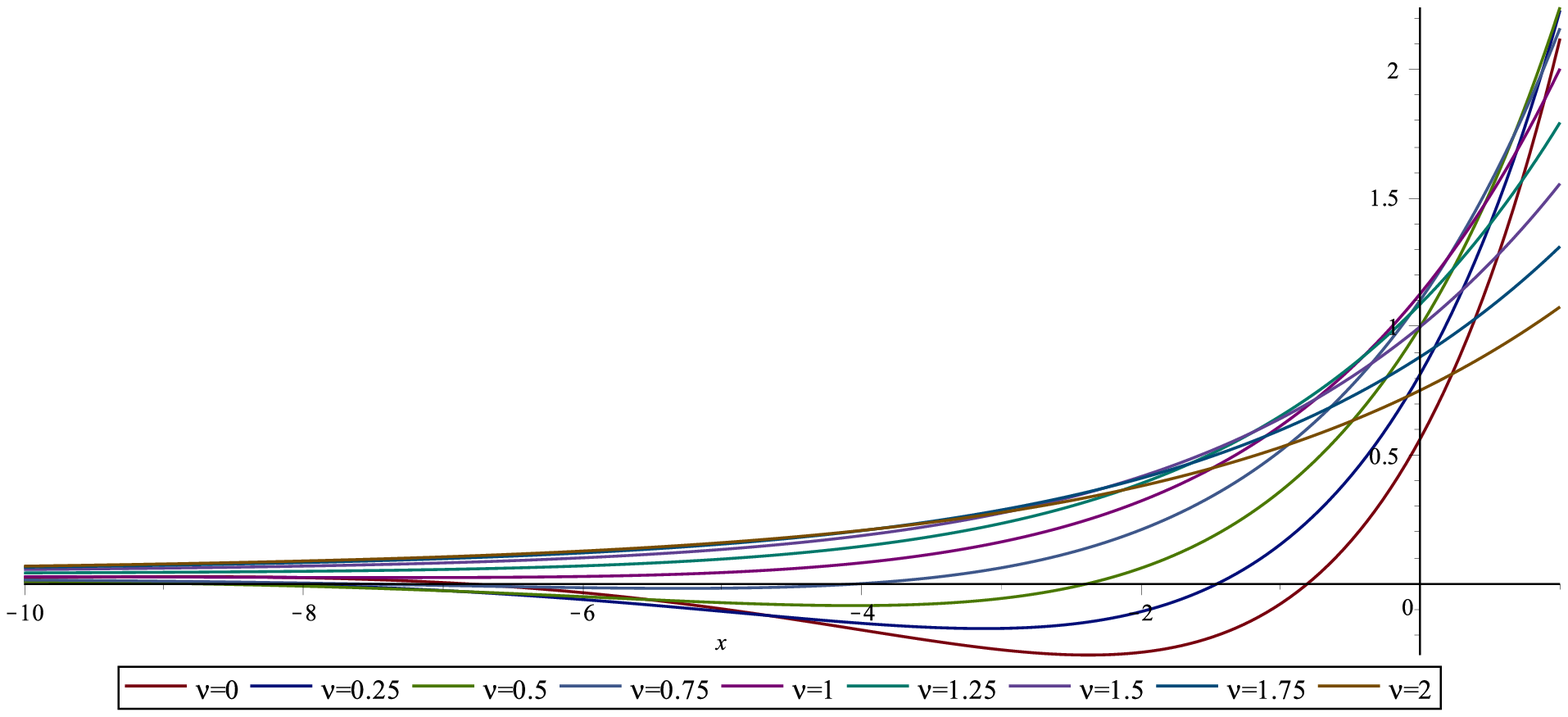}}}
\caption{}
\label{fig1}
\end{figure}

In Addition, we present some examples of the 3-parameters function $\mathcal{W}_{\alpha,\beta,\nu}$ in the following table, and in Figure \ref{fig1} we represent the behavior of this function for different values of the parameters $\alpha,\beta,\nu$;

\begin{center}
  \begin{tabular}{ | l | c |  }
    \hline
    Integer order derivatives & Fractional order derivatives  \\ \hline
    $\mathcal{W}_{0,1,0}(x)=e^x$ & $\mathcal{W}_{\frac{1}{2},\frac{1}{2},\frac{1}{2}}(\sqrt x)=+I_{0}(2\sqrt x)+L_{0}(2\sqrt x)$ \\ \hline
    $\mathcal{W}_{0,1,n}(x)=\frac{e^x}{x^n}-\sum_{i=0}^{n-1}\frac{x^{i-n}}{i!}$ with $n\in\mathbb{N}$& $\mathcal{W}_{\frac{1}{2},\frac{1}{2},\frac{3}{2}}(\sqrt x)=+I_{1}(2\sqrt x)+L_{1}(\sqrt x)$  \\ \hline
    $\mathcal{W}_{1,1,0}(x)=\sqrt x I_{1}(2\sqrt x)$& $\mathcal{W}_{\frac{1}{2},\frac{1}{2},1}(\sqrt x)=\frac{\sinh(2\sqrt x)+\cosh(2\sqrt x)-1}{\sqrt{\pi x}}$ \\ \hline
$\mathcal{W}_{1,1,\nu}(x)=x^{-\frac{\nu-1}{2}}I_{\nu-1}(2\sqrt x)$ & $\mathcal{W}_{\frac{1}{2},\frac{1}{2},2}(\sqrt x)=\frac{\left(2\sqrt x-1\right) e^{2\sqrt x}-2x+1}{2x\sqrt{\pi x}}$ \\ \hline
  \end{tabular}
\end{center}

where $I_\alpha(x)=i^{-\alpha}J_\alpha(ix)=\sum_{m=0}^\infty\frac{1}{m!\Gamma(m+\alpha+1)}(\frac{x}{2})^{2m+\alpha}$ is the modified Bessel function of the first kind and $L_\alpha(x)=\left(\frac{x}{2}\right)^{\nu+1}\sum_{m=0}^\infty\frac{\left(\frac{x}{2}\right)^{2m}}{\Gamma(m+\frac{3}{2})\Gamma(m+\nu+\frac{3}{2})}$ is the modified Struve function.

\section{Laplace Transform}
Let us compute the Laplace transform of the $\mathcal{W}^{\left(\bar{\alpha},\bar{\nu}\right)}(\lambda x)$

\begin{eqnarray}\label{LT}
\mathcal{L}\left(\mathcal{W}^{\left(\bar{\alpha},\bar{\nu}\right)}(\lambda x^{\alpha_{n+1}}),s\right)&=&\int_{0}^{\infty} e^{-sx} \sum_{k=0}^{\infty}\prod_{i=1}^{k} \prod_{j=1}^{n} \frac{\Gamma(\alpha_{n+1}i+a_j)}{\Gamma(\alpha_{n+1}i+b_j)}\frac{\lambda^k x^{\alpha_{n+1}k}}{\Gamma(\alpha_{n+1}k+b_{n+1})}dx\nonumber\\
&=&\sum_{k=0}^{\infty}\prod_{i=1}^{k} \prod_{j=1}^{n} \frac{\lambda^k \Gamma(\alpha_{n+1}i+a_j)}{\Gamma(\alpha_{n+1}i+b_j)\Gamma(\alpha_{n+1}k+b_{n+1})}\int_{0}^{\infty} e^{-sx}x^{\alpha_{n+1}k} dx\nonumber\\
&=&\frac{1}{s}\sum_{k=0}^{\infty}\prod_{i=1}^{k} \prod_{j=1}^{n} \frac{\Gamma(\alpha_{n+1}i+a_j) \Gamma(\alpha_{n+1}k+1)}{\Gamma(\alpha_{n+1}i+b_j)\Gamma(\alpha_{n+1}k+b_{n+1})}\left(\frac{\lambda}{s^{\alpha_{n+1}}}\right)^k.
\end{eqnarray}
th analytical properties of the $\mathcal{W}^{\left(\bar{\alpha},\bar{\nu}\right)}$ provides that the resulting Laplace transform turns out to be an analytic function, vanishing at infinity and exhibiting an essential singularity at $s=0$. 
\begin{remark}
In case we set $\alpha_j=1,\,\,\,j=1,...,n+1$, the multi-index special functions $\mathcal{W}^{\left(\bar{\alpha},\bar{\nu}\right)}$ will be related to the hyper-Bessel functions as is showed in (\ref{HB}). After some calculations, we obtain the following functional relation between the Laplace transform of the Hyper-Bessel function and the multi-index Mittag-Leffler function. A more general relation between these two functions can be found in the article of Kiryakova and Luchko \cite{KirLuch10}.

\begin{eqnarray}\label{LTHB}
\mathcal{L}\left(\mathcal{W}^{\left(\bar{1},\bar{\nu}\right)}(\lambda x),s\right)&=&\prod_{j=1}^n\Gamma(1+a_j)\frac{1}{s}\sum_{k=0}^{\infty}\frac{1}{\prod_{j=1}^n \Gamma(k+a_{j+1}+1)}\left(\frac{\lambda}{s}\right)^k\nonumber\\
&=&\prod_{j=1}^n\Gamma(1+a_j)\frac{1}{s}E_{(1,1,...,1),(a_{j+1}+1)}^{(n)}\left(\frac{\lambda}{s}\right).
\end{eqnarray}

\end{remark}

\begin{remark}
The Laplace transform of $\mathcal{W}_{\alpha,\beta,\nu}(x)$ can be obtained as a special case of the (\ref{LT}) as follows:

\begin{equation}
\mathcal{L}\left(	\mathcal{W}_{\alpha,\beta,\nu}(\lambda x^\rho),s\right)	=\frac{1}{s}\sum_{k=0}^\infty\prod_{i=1}^k\frac{\beta\Gamma(\beta i+1-\alpha)}{\Gamma(\beta i+1)}\frac{\Gamma(\rho k+1)}{\Gamma(\beta k+1-\alpha+\nu)}\left(\frac{\lambda}{s^\rho}\right)^k;
\end{equation}

and, in the case of the parameter $\alpha=1$, we obtain the well-known Laplace transform of the Wright function which can be expressed in terms of the two-parameter Mittag-Leffler function.

\begin{equation}
\mathcal{L}\left(	\mathcal{W}_{1,\beta,\nu}(\lambda x),s\right)=\mathcal{L}\left(	W_{\beta,\nu}(\lambda x),s\right)=\frac{1}{s}E_{\beta,\nu}\left(\frac{\lambda}{\beta\,s}\right).
\end{equation}•
\end{remark}

\section{Recurrence relations of $\mathcal{W}_{\alpha,\beta,\nu}$}

A recurrence relation is an equation that recursively defines a sequence of values; given one or more initial terms, each further term of the sequence is defined as a function of the previous terms. Differential recurrence relation of the generalized Wright function can be used in the study of fractional differential equations, and it is obtained directly from series representation.  

\begin{equation}\label{RR1}
x^{\alpha+\beta}\frac{d^\beta}{dx^\beta}\left(x^{\nu-\alpha+\beta}\mathcal{W}_{\alpha,\beta,\nu+\beta}(x^\beta)\right)-2x^{\nu+\beta}\mathcal{W}_{\alpha,\beta,\nu}(x^\beta)+x^{\alpha+\nu}\frac{d^\alpha}{dx^\alpha}\mathcal{W}_{\alpha,\beta,\nu-\beta}(x^\beta)=0.
\end{equation}

\begin{remark}\label{RRB-C}
In case $\alpha=\beta=1$; $\nu=n+1$  with $n\in\mathbb{N}_0$ and 
$$\frac{d}{dx}\mathcal{W}_{1,1,n}(x)=\mathcal{W}_{1,1,n+1}(x)=C_n(x);$$ 
we obtain the well known three-term recurrence relation for the \textit{Bessel-Clifford function} $C_n(x)$
\begin{equation}\label{RRBC}
xC_{n+2}(x)+(n+1)C_{n+1}(x)=C_n(x).
\end{equation}
\end{remark}

\begin{remark}{\textbf{Recurrence fractional derivatives relation for the Wright and Mittag-Leffler functions.}}
From the relation (\ref{Wright}) between the generalized Wright function and the classical Wright function the relation (\ref{RR1})  becomes

\begin{equation}
 \frac{d^\lambda}{dz^\lambda}\left(z^{\lambda+\nu-1}W_{\lambda,\lambda+\nu}\left(\frac{z^\lambda}{\lambda}\right)\right)=z^{\nu-1}W_{\lambda,\nu}\left(\frac{z^\lambda}{\lambda}\right);
\end{equation}

using the formula
\begin{equation}
 \frac{d}{dz}W_{\lambda,\nu-\lambda}\left(\frac{z^\lambda}{\lambda}\right)=z^{\lambda-1}W_{\lambda,\nu}\left(\frac{z^\lambda}{\lambda}\right).
\end{equation}

In case $\alpha=0$, and $\beta\rightarrow\alpha,\,\,\nu\rightarrow\beta-1$ the generalized Wright function is related to the Mittag-Leffler function by the following relation:
\begin{equation}
\mathcal{W}_{0,\alpha,\beta-1}(z)=E_{\alpha,\beta}(z).
\end{equation}
In particular, from the recurrence relation (\ref{RR1}), we obtain the new recurrence relation involving fractional derivatives for M-L functions.

\begin{equation}
z^\alpha \frac{d^\alpha}{dz^\alpha}\left(z^{\alpha+\beta-1}E_{\alpha,\alpha+\beta}(z^\alpha)\right)-2z^{\alpha+\beta-1}E_{\alpha,\beta}(z^\alpha)+z^{\beta-1}E_{\alpha,\beta-\alpha}(z^\alpha)=0.
\end{equation}
\end{remark}

\section{Partial derivatives of $\mathcal{W}_{\alpha,\beta,\nu}$ with respect to the parameters}
In this section, taking inspiration from the works of Apelblat and Mainardi \cite{A20}, \cite{AM20} we analyse the derivatives of $\mathcal{W}_{\alpha,\beta,\nu}$ respect the three parameters included in the function. We can treat parameters as variables and hence the derivatives with respect to them can be obtained. These derivatives lead to infinite power series involving digamma ($\psi$) and gamma functions.

\begin{equation}
\frac{\partial}{\partial\nu}\mathcal{W}_{\alpha,\beta,\nu}(z)=-\sum_{k=0}^\infty\prod_{i=1}^k\frac{\Gamma(\beta i+1-\alpha)}{\Gamma(\beta i+1)}\frac{\psi(\beta k+1-\alpha+\nu)}{\Gamma(\beta k+1-\alpha+\nu)}z^k;
\end{equation}

\begin{eqnarray}
&&\frac{\partial}{\partial\beta}\mathcal{W}_{\alpha,\beta,\nu}(z)=\sum_{k=0}^\infty\prod_{i=1}^k\frac{\Gamma(\beta i+1-\alpha)}{\Gamma(\beta i+1)}\frac{z^k}{\Gamma(\beta k+1-\alpha+\nu)}\cdot \nonumber\\
&&\cdot\left[\sum_{j=1}^kj\left[\psi(\beta j+1-\alpha)-\psi(\beta j+1)\right]-k\psi(\beta k+1-\alpha+\nu)\right];
\end{eqnarray}

\begin{equation}
\frac{\partial}{\partial\alpha}\mathcal{W}_{\alpha,\beta,\nu}(z)=\sum_{k=0}^\infty\prod_{i=1}^k\frac{\Gamma(\beta i+1-\alpha)}{\Gamma(\beta i+1)}\frac{z^k}{\Gamma(\beta k+1-\alpha+\nu)}\left[-\sum_{j=1}^k\psi(\beta j+1-\alpha)+\psi(\beta k+1-\alpha+\nu)\right];
\end{equation}

where $\psi(z)=\frac{\Gamma'(z)}{\Gamma(z)}$ denotes the \textit{digamma} function.

\begin{remark}
In the case $\alpha=1$ and considering the property of the digamma function $\psi(z+1)=\psi(z)+\frac{1}{z}$; we obtain
the formula (5) and (6) of the Apelblat-Mainardi article (\cite{AM20}) for the classical Wright function
$$\frac{\partial}{\partial\beta}\mathcal{W}_{1,\beta,\nu}(z)=\left(\frac{\partial}{\partial\beta}W_{\beta,\nu}\right)(\beta z)=-\sum_{k=0}^\infty\left(\frac{\psi(\beta k+\nu)}{k!\Gamma(\beta k+\nu)}\right)kz^k;$$ 
$$\frac{\partial}{\partial\nu}\mathcal{W}_{1,\beta,\nu}(z)=\left(\frac{\partial}{\partial\nu}W_{\beta,\nu}\right)(\beta z)=-\sum_{k=0}^\infty\left(\frac{\psi(\beta k+\nu)}{k!\Gamma(\beta k+\nu)}\right)z^k.$$ 

By setting the parameters, $\alpha=0,\,\,\beta\rightarrow\alpha$ and $\nu\rightarrow\beta-1$, we obtain the formulas (95) and (96) of the Apelblat paper (\cite{A20}) 
$$\frac{\partial}{\partial\alpha}\mathcal{W}_{0,\alpha,\beta-1}(z)=\frac{\partial}{\partial\alpha}E_{\alpha,\beta}(z)=-\sum_{k=0}^\infty\left(\frac{k\psi(\alpha k+\beta)}{\Gamma(\alpha k+\beta)}\right)z^k;$$ 
$$\frac{\partial}{\partial\beta}\mathcal{W}_{0,\alpha,\beta-1}(z)=\frac{\partial}{\partial\beta}E_{\alpha,\beta}( z)=-\sum_{k=0}^\infty\left(\frac{\psi(\alpha k+\beta)}{\Gamma(\alpha k+\beta)}\right)z^k.$$ 
\end{remark}

\section{Conclusion}

The aim of this paper is  to investigate several properties related to the  multi-index special function $\mathcal{W}^{\left(\bar{\alpha},\bar{\nu}\right)}$ and its 3-parameters version. An important result was finding the connection between the $\mathcal{W}^{\left(\bar{\alpha},\bar{\nu}\right)}$ and the  hyper-Bessel function of Delerue.
Here we analyzed the Laplace transform, recurrence relation and derivatives of the function with respect to the parameters. In particular, we found new findings that, for special values of the parameters, retrieve some well-known relations. Indeed, a simple functional relation is obtained between the Laplace transform of the hyper-Bessel function and the multi-index Mittag-Leffler. 
\section*{Disclosure statement}
No potential conflict of interest was reported by the author.

\section*{Acknowledgements}
The author is grateful to Dr Roberto Garra for providing essential
information, help and advice.


\appendix


\section{Fractional calculus}

In order to make the  papar self-contained, we briefly recall main definitions and properties of fractional calculus operators.

Let $\gamma\in \mathbb{R}^{+}$. The Riemann-Liouville fractional integral is defined by

\begin{equation}
J^{\gamma}_x f(x) =
\frac{1}{\Gamma(\gamma)}\int_0^{x}(x-x')^{\gamma-1}f(x') dx',
\label{riemann-l}
\end{equation}
where
$$\Gamma(\gamma)= \int_0^{+\infty}x^{\gamma-1}e^{-x}dx,$$
is the Euler Gamma function.\\

Note that, by definition, $J^0_x f(x)= f(x)$. \\

Moreover it satisfies the semigroup property, i.e. $J_x^{\alpha}J_x^{\beta} f(x)= J_x^{\alpha+\beta}f(x)$.\\
There are different definitions of fractional derivative (see e.g. \cite{Pod99}). In this paper we used the fractional derivatives in the sense of Caputo, that is

\begin{equation}
D_x^{\gamma}f(x)=   J_x^{m-\gamma} D_x^m f(x)=
\frac{1}{\Gamma(m-\gamma)}\int_0^{x}(x-x')^{m-\gamma-1}\frac{d^m}{d(x')^m}f
(x') \, \mathrm dx', \;\gamma \ne m.
\end{equation}

It is simple to prove the following properties of fractional
derivatives and integrals (see e.g. \cite{Pod99}) that will be used
in the analysis:
\begin{align}
&D_x^{\gamma} J_x^{\gamma} f(x)= f(x), \quad \gamma> 0,\\
&J_x^{\gamma} D_x^{\gamma} f(x)= f(x)-\sum_{k=0}^{m-1}f^{(k)}(0)\frac{x^k}{k!}, \qquad \gamma>0, \: x>0,\\
&J_x^{\gamma} x^{\delta}= \frac{\Gamma(\delta+1)}{\Gamma(\delta+\gamma+1)}x^{\delta+\gamma} \qquad \gamma>0, \: \delta>-1, \: t>0,\\
&D_x^{\gamma} x^{\delta}=
\frac{\Gamma(\delta+1)}{\Gamma(\delta-\gamma+1)}x^{\delta-\gamma}
\qquad \gamma>0, \: \delta>-1, \: t>0.
\end{align}

\appendix

\end{document}